\documentclass[]{elsarticle}

\RequirePackage{amsthm,amsmath,amssymb,bbm}
\RequirePackage[numbers]{natbib}

 \usepackage{rotating}
\usepackage{colortbl}
\usepackage{multirow}
\usepackage[table]{xcolor}
\usepackage{makecell}

\numberwithin{equation}{section}
\theoremstyle{plain}
\newtheorem{theorem}{Theorem}
\newtheorem{assume}{Assumption}
\newtheorem{lemma}{Lemma}
\newtheorem{corollary}{Corollary}
\newtheorem{prop}{Proposition}
\newtheorem{remark}{Remark}

\usepackage{orcidlink}

\allowdisplaybreaks

\begin{document}
\begin{frontmatter}
	
\title{PAC-Bayesian risk bounds for fully connected deep neural network with Gaussian priors}
\author[aaatienmt]{The Tien Mai\,\orcidlink{0000-0002-3514-9636} }\ead{the.tien.mai@fhi.no}

\affiliation[aaatienmt]{organization={
Norwegian Institute of Public Health},
city={Oslo},
postcode={0456}, 
country={Norway}}

\begin{abstract}
Deep neural networks (DNNs) have emerged as a powerful methodology with significant practical successes in fields such as computer vision and natural language processing.  Recent works have demonstrated that sparsely connected DNNs with carefully designed architectures can achieve minimax estimation rates under classical smoothness assumptions. However, subsequent studies revealed that simple fully connected DNNs can achieve comparable convergence rates, challenging the necessity of sparsity. Theoretical advances in Bayesian neural networks (BNNs) have been more fragmented. Much of those work has concentrated on sparse networks, leaving the theoretical properties of fully connected BNNs underexplored. 
In this paper, we address this gap by investigating fully connected Bayesian DNNs with Gaussian prior using PAC-Bayes bounds. We establish upper bounds on the prediction risk for a probabilistic deep neural network method, showing that these bounds match (up to logarithmic factors) the minimax-optimal rates in Besov space, for both nonparametric regression and binary classification with logistic loss. Importantly, our results hold for a broad class of practical activation functions that are Lipschitz continuous. 
\end{abstract}


\end{frontmatter}

\section{Introduction}
\label{sc_introduction}
Deep neural networks have become one of the most powerful and versatile tools in modern machine learning, with widespread applications in pattern recognition and nonparametric regression \citep{anthony1999neural,haykin1998neural,hertz2018introduction,ripley2007pattern}. Over the past decade, the advent of deep learning—driven by multilayer feedforward neural networks with numerous hidden layers—has enabled groundbreaking achievements in domains such as computer vision and natural language processing \citep{goodfellow2016deep,lecun2015deep}. These practical successes have sparked a surge of interest in the theoretical analysis of neural networks, with particular focus on their generalization properties and approximation capabilities. Foundational works, including \cite{Barron94estimation,Bartlett2017MarginBoundsNNs,Srebro2018PACMarginBoundsNNs,ZhangUnderstandingDL2017,SchmidtHieberDNN,Suzuki18DNNkerels,Imaizumi19DNN,suzuki2019adaptivity,Suzuki2019Superiority,elbrachter2021deep}, have made significant strides in analyzing these properties, laying the groundwork for rigorous study.  

In a seminal contribution, \cite{SchmidtHieberDNN} demonstrated that sparsely connected deep neural networks (DNNs) with ReLU activation functions and carefully constructed architectures can achieve minimax-optimal estimation rates (up to logarithmic factors) under classical smoothness conditions. This breakthrough stimulated a surge of interest in the theoretical analysis of Bayesian neural networks (BNNs) \cite{Rockova2018}, with the majority of subsequent work focusing on sparse architectures—where sparsity typically refers to the presence of a limited number of non-zero weights. For some time, such sparsity was viewed as essential for attaining desirable convergence rates in both frequentist and Bayesian frameworks.

However, recent developments have begun to challenge this prevailing assumption. Notably, a series of studies including \cite{kohler2021rate}, \cite{langer2021analysis}, and \cite{kohler2023estimation} have shown that fully connected DNNs, without explicit sparsity constraints, can also achieve comparable minimax rates. These findings open up new theoretical avenues and encourage further exploration of the performance and generalization properties of fully connected neural networks, for example \cite{yara2024nonparametric}. In particular, they motivate continued investigation within Bayesian frameworks.

On the Bayesian side, theoretical research on deep neural networks has been comparatively limited and fragmented. Several studies have explored the properties of posterior distributions in Bayesian neural networks, including \cite{Rockova2018}, \cite{Suzuki18DNNkerels}, \cite{liu2021variable,sun2022learning}, \cite{lee2022asymptotic}, and \cite{kong2023masked}. In a related line of work, \cite{Vladimirova2019PriorsBNNsUnits} investigated the regularization effects induced by priors at the unit level. Additional contributions, such as \cite{bai2020efficient,jantre2023layer,jantre2024spike} and \cite{bhattacharya2024comprehensive}, have analyzed aspects of variational inference in Bayesian neural networks. However, much of this literature remains centered on sparse network architectures, with fully connected Bayesian neural networks largely unexplored from a theoretical perspective. A recent exception is the work by \cite{kong2024posterior}, which examined the frequentist properties of the posterior distribution in fully connected neural networks, thereby underscoring the importance and timeliness of further research in this direction. Importantly, while these studies provide valuable insights into the frequentist behavior of the posterior, they do not address the predictive performance of Bayesian deep neural networks. Our work aims to help fill this gap by adopting a PAC-Bayesian approach to quantify prediction risk in fully connected architectures.

In this paper, we derive upper bounds on the prediction risk for an exponentially weighted aggregate estimator and demonstrate that these bounds achieve, up to logarithmic factors, the minimax-optimal convergence rates in both nonparametric regression and binary classification with logistic loss. Our theoretical results are established under the widely used Gaussian prior \cite{bingham2019pyro}, thereby providing a rigorous foundation for the analysis of fully connected Bayesian neural networks. The PAC-Bayesian inequality framework, originally introduced to obtain generalization bounds for learning algorithms \cite{McA, STW} (where PAC stands for “Probably Approximately Correct”), has since evolved into a versatile tool for statistical learning theory. Notably, \cite{catonibook} extended its application to the derivation of oracle inequalities and excess risk bounds, a perspective we adopt in our analysis. The practical relevance of PAC-Bayesian bounds has been demonstrated in a variety of recent studies, including \cite{mai2023reduced, mai2024high, mai2024sparse}. For comprehensive introductions and recent developments in PAC-Bayesian theory, we refer the reader to \cite{alquier2021user} and \cite{guedj2019primer}.

\textit{ More on related works and our contributions}: 
Previous research on PAC-Bayesian bounds for neural networks has predominantly concentrated on sparse architectures. For example, \cite{cherief2020convergence} analyzed PAC-Bayesian bounds for variational approximations in sparse Bayesian neural networks within nonparametric regression frameworks; while \cite{steffen2022pac} investigated similar bounds using Gibbs posteriors, also in the context of regression. In the classification setting, \cite{mai2024misclassification} derived misclassification error bounds for sparse Bayesian neural networks based on hinge loss. More recently, \cite{alquier2024minimaxopti} established PAC-Bayesian bounds for sparse neural networks under dependent data structures. In contrast to these studies, the present work focuses on fully connected deep neural networks, thereby broadening the applicability of PAC-Bayesian theory beyond sparsity-constrained models. In a related direction, \cite{tinsi2022risk} provided risk bounds for shallow neural networks—specifically, networks with a single hidden layer—equipped with Gaussian priors, in regression settings. Our contribution extends these results to deeper, fully connected architectures and, importantly, incorporates both regression and binary classification scenarios.

\paragraph{Structure of the paper} Section \ref{sc_nota_defini} presents the fundamental notations and terminology that will be used throughout the paper. In this section, we also provide a formal definition of the fully connected deep neural network architecture that forms the basis of our theoretical development. In Section \ref{sc_regression}, we introduce our proposed methodology and establish risk bounds in the context of nonparametric regression, highlighting the theoretical guarantees associated with the performance of the estimator. 
Section \ref{sc_classification} extends the analysis to the binary classification setting. Here, we derive bounds on the prediction error for the cross-entropy loss and, in particular, present a novel result concerning the misclassification error, which contributes to the theoretical understanding of classification performance in deep learning models. Finally, \ref{sc_proof} contains all technical proofs.

\section{Notations and definitions}
\label{sc_nota_defini}
\subsection{Notations}

Here, we introduce notations used in this paper.  For any vector $x=(x_1,...,x_d) \in \Omega := [0,1]^d$ and any real-valued function $f$ defined on $\Omega $, we denote:
$
\|x\|_\infty = \max_{1\leq i \leq d} |x_i|
\textnormal{,}
\,
\textnormal{ and }
\,
\|f\|_\infty = \sup_{t\in \Omega } |f(t)|.
$
For any $\mathbf{k}\in \{0,1,2,...\}^d$, we define $|\mathbf{k}|=\sum_{i=1}^d k_i$ and the mixed partial derivatives when all partial derivatives up to order $|\mathbf{k}|$ exist:
$
D^{\mathbf{k}} f(x) = \frac{\partial^{|\mathbf{k}|}f}{\partial^{k_1}x_1...\partial^{k_d}x_d} (x) .
$ Let $ \beta > 0$ be the smoothness parameter and $m = \lfloor \beta \rfloor$. The $\beta $-H\"older norm of function $f$ is defined by
$$
\| f \|_{C^\beta} 
:= 
\sup_{\| u\|_1 \leq m} \| D^u f\|_{\infty} + \sup_{\| u\|_1 = m} \sup_{x, y \in \Omega} \frac{|D^u f(x)  - D^u f(y)|}{|x - y|^{\beta - m}}.
$$
The $\beta$-H\"older space $ \mathcal{H}^\beta(\Omega)$ is defined as a set of functions with bounded $\beta$-H\"older norm. 

Let $k \in \mathbb{N}$ and $1 \leq p \leq \infty$. The Sobolev norm of a function $f$ is defined by 
$$
\| f \|_{W^{k, p}} := 
\big( \sum_{\| \alpha\|_1 \leq k} \| D^\alpha f \|_{L^p}^p \big)^{1/p}, 1 \leq p < \infty 
$$ 
and
$ \| f \|_{W^{k, \infty}} = \sup_{\| \alpha\|_1 \leq k} \| D^\alpha f\|_{\infty} $. 
The Sobolev space $W^{k, p}(\Omega)$ is defined by the function space consisting of functions with bounded Sobolev norm.

Note that the notion of the smoothness of a function is related to the differentiability of the function. The Besov space extends the concept of smoothness.
Before defining the Besov space, we define the $r$-th modulus of smoothness of the function $f$,  \citep{gine2016mathematical}, as: 
$
w_{r, p}(f, t) = \sup_{ \|h\|_2 \leq t} \| \Delta_h^r (f)\|_p,
$ where
\begin{equation*}
	\Delta_h^r(f)(x) 
	= 
	\begin{cases} 
		\sum_{j=0}^r \binom{r}{j} (-1)^{r-j}f(x+jh) & x \in \Omega, x+ rh \in \Omega, 
		\\
		0 & \text{otherwise } 
		.
	\end{cases}
\end{equation*} 
For $0 < p, q \leq \infty,~\beta >0,~r=\lfloor \beta \rfloor +1$, define the Besov norm of a function $f$ by
\begin{equation*}
	\| f \|_{B_{p,q}^\beta}
	:= 
	\|f\|_p + \begin{cases} \left(\int_0^\infty \left( t^{-\beta}w_{r,p}(f,t)\right)^q \frac{dt}{t}\right)^{1/q} & q < \infty, 
		\\ 
		\sup_{t>0} \{t^{-\beta}w_{r,p}(f,t)\} & q = \infty.  \end{cases}
\end{equation*}
The Besov space $B_{p,q}^\beta(\Omega)$ is defined as the set of functions with finite Besov norms \citep{gine2016mathematical}, i.e., 
$
B_{p,q}^\beta (\Omega) = \{f: \| f\|_{B_{p,q}^\beta } < \infty \}.
$

Note that the Besov spaces can be defined without the differentiability and continuity of functions, and are more general than the H\"{o}lder and Sobolev spaces, which are subspaces of the Besov spaces. The Besov space may contain complicated functions including discontinuous function if $d/p \geq \beta $. For example, the Cantor function does not belong to any Sobolev space since it does not have a weak derivative but belongs to a Besov space \citep{sawano2018theory}.

Let $ R_1,R_2 $ be two probability measures. The  Kullback-Leibler divergence is defined by
$
\mathcal{K}(R_1,R_2)  = 
\int \log \left(\frac{{\rm d} R_1 }{{\rm d} R_2} \right){\rm d} R_1 $  if  $ R_1 \ll R_2
$, and  $
+ \infty \text{ otherwise}.
$

\subsection{Fully connected deep neural networks}

A fully connected deep neural network can be described as a function \( f_\theta: \mathbb{R}^d \rightarrow \mathbb{R} \), defined recursively in the following manner that with $ x^{(0)}  = x $, put
\begin{align*}
	f_\theta(x) & = \phi( x^{(L-1)}) ,
	\\
	x^{(\ell)} & = \phi(W_\ell x^{( \ell-1)} + b_\ell) \quad \text{for} \quad \ell = 1, \ldots, L-1 ,
\end{align*}
where:
\begin{itemize}
	\item \( L\geq 3 \) is the number of layers,
	\item \( \theta = \{W_\ell, b_\ell\}_{\ell=1}^L \) are the parameters of the network, with \( W_\ell \in \mathbb{R}^{D_\ell \times D_{\ell-1}} \) being the weight matrix and \( b_\ell \in \mathbb{R}^{D_\ell } \) being the bias vector for layer \( \ell \),
	\item \( \phi(\cdot) \) is an activation function (e.g., ReLU, sigmoid, tanh) acting componentwise,
	\item \( x^{(\ell)} \) is the output (hidden state) of the \( \ell \)-th layer.
\end{itemize}

In this recursive definition, each layer transforms its input \( x^{(\ell-1)} \) through a linear transformation followed by a non-linear activation function $\phi $, culminating in the final output \( f_\theta(x) \). Each weight matrix \( W_\ell \in \mathbb{R}^{D_\ell \times D_{\ell-1}} \) contains entries \( (i,j) \) that are edge weights, connecting the \( j \)-th neuron in layer \( \ell-1 \) to the \( i \)-th neuron in layer \( \ell \). The vector \( b_\ell \in \mathbb{R}^{D_\ell} \), known as the bias vector, has entries where the \( i \)-th element corresponds to the bias for the \( i \)-th neuron in layer \( \ell \).

Let \( D_0 = d \) represent the number of units in the input layer, \( D_L = 1 \) the number of units in the output layer, and \( D_\ell = D \) the number of units in the hidden layers. The network architecture is defined by its depth, \( L \geq 3 \), and its width, \( D \geq 1 \). For a fully connected network, the total number of coefficients is given by \( T := LD(D+1) \).  

We assume that \( L \) and \( D \) are fixed and require that \( d \leq D \). 
The parameter for a fully connected Deep Neural Network (DNN) is denoted as \( \theta = \{ (W_1, b_1), \ldots, (W_L, b_L) \} \), and we define \( \Theta_{L,D} \) as the set of all possible parameters. Alternatively, we will also consider the stacked coefficients parameter, \( \theta = (\theta_1, \ldots, \theta_T) \).

We use in this paper independent Gaussian priors $ \mathcal{N}(0,1) $ for each parameter's component. However, the results remain valid for any Gaussian prior of the form $\mathcal{N}(0, \sigma^2) $ with known variance $\sigma^2 $.

In this paper, our results are shown to hold for a wide variety of activation functions, covering popular options like the ReLU, sigmoid, tanh activation and the identity map. This generality highlights the robustness of our approach to different neural network designs. To proceed with the theoretical development, we now specify the following assumption related to the fully connected deep neural networks.

\begin{assume}
	\label{asm1}
	Assume that the activation function $\phi$ is $1$-Lispchitz  with respect to the absolute value: 
	$$
	|\phi(x)|\leq|x|
	$$ 
	for any $x\in\mathbb{R}$.  Moreover, for $ C\ge1$ we assume that 
	$\Vert f\Vert_{\infty}\le C$. It is assumed that the absolute values of all coefficients are bounded above by a constant \( B \geq 1 \).
\end{assume}

It is noted that the boundedness assumptions above is standard as it has been used for theoretical studies of deep neural networks in various works such as \cite{Rockova2018,SchmidtHieberDNN,kim2021fast}.

\section{Nonparametric regression}
\label{sc_regression}
\subsection{Problem setting}
We consider the nonparametric regression framework. We have a collection of random samples $( Y_i, X_i) \in \mathbb{R} \times [0,1]^d  $ for $i=1,...,n$ which are independent and identically distributed (i.i.d.) copies of generic random variables $ (Y, X) $
on some probability space $(\Omega,\mathcal{A},\mathbb{P} ) $ with the following relationship:
$$
Y = f_0(X)+ \epsilon ,
$$
where $ X $ follows the uniform distribution on the set $[0,1]^d $,
and
$\epsilon $ is random noise variable with mean $\mathbb{E}\, (\epsilon \, | X ) =0 $. 
Here, 
$f_0:[0,1]^d\rightarrow \mathbb{R}$ is the true unknown function. For instance, the true regression function $f_0$ may belong to the set $ \mathcal{H}^\beta $ of H\"older functions with level of smoothness $\beta$.

For any DNN estimator $ f_\theta $, the prediction risk and its empirical counterpart
are given by 
$$
R( \theta ) 
:= 
R( f_\theta ) 
:= 
\mathbb{E} \big[(Y- f_\theta (X))^{2}\big]
\quad\text{and}
\quad 
r_{n}( \theta )
=
r_{n}( f_\theta )
=
\frac{1}{n}\sum_{i=1}^{n}
\big(Y_{i} - f_\theta (X_{i})\big)^{2}
,
$$
respectively, where $\mathbb{E} $ denotes the expectation under $ \mathbb{P} $ and $\mathbb{E}_Z$ is the (conditional) expectation only with respect to a random variable $Z$. We work with random design that $ X  \sim  \mathrm{P}_{ X }, $ where $\mathrm{P}_{ X }$ is the marginal distribution of $X$.  The accuracy of the estimation procedure will be quantified
in terms of the excess risk 
\begin{equation*}
	R( f_\theta )-R(f_0)
	=
	\mathbb{E}_{X} \big[(f_\theta (X)-f_0(X))^{2}\big]
	=
	\| f_\theta - f_0 \|_{L^{2}( \mathrm{P}_{X})}^{2}
	.
\end{equation*}

We adopt a PAC-Bayesian approach \citep{alquier2021user} and consider an exponentially weighted aggregate (EWA) procedure. Let's consider the following Gibbs posterior distribution:
\begin{align}
	\label{eq_mainporsterior}
	\hat{\rho}_\lambda(\theta)
	\propto
	\exp[-\lambda r_n(\theta)] \pi(\theta)
\end{align}
where $\lambda>0$ is a tuning parameter that will be discussed later and $\pi(\theta)$ is a prior distribution on the parameter  $ \theta $. It is worth noting that the Gibbs posterior defined in \eqref{eq_mainporsterior} inherently favors parameters with low empirical risk \( r_n \), due to the exponential weighting of the negative risk. As a result, parameters yielding smaller empirical risk values are assigned higher posterior probability mass. 

When the noise \( \epsilon \) follows a Gaussian distribution \( \mathcal{N}(0, \sigma^2) \), setting \( \lambda = n / 2\sigma^2 \) in \eqref{eq_mainporsterior} yields the standard Bayesian posterior distribution. However, our approach—through the use of a general tuning parameter \( \lambda \)—extends naturally to broader noise settings, such as those described in Assumption \ref{assu_bounded}. This flexibility allows our framework to accommodate non-Gaussian noise distributions, thereby positioning it within the broader class of quasi-Bayesian methods.

\subsection{Minimax-optimal excess risk bound in the Besov space}

Let $\theta^*$ denote a minimizer of $R$ when it exists:
$
R(\theta^*) 
= 
\min_{\theta\in \Theta_{L,D} } R(\theta). 
$

\begin{assume}
	\label{assu_bounded}
	Assume that there exist two constants $\sigma, \varsigma > 0 $ such that 
	$$
	\mathbb{E} ( \, \vert \epsilon \vert^{k} | X )
	\le
	\frac{k!}{2}\sigma^{2} \varsigma^{k-2}
	\quad 
	\text{  for all $ k\ge 2 $ }
	.
	$$
\end{assume}

From Assumption \ref{assu_bounded}, it is showed that the probability tails of the random noise \( \epsilon \) decay as fast as those of a Gamma distribution. This kind of assumption is quite common in machine learning, as it ensures the noise does not exhibit heavy tails. Importantly, Gaussian variables satisfy this assumption, as do centered random variables with values bounded by a fixed constant (see Chapter 2 of \cite{boucheron2013concentration}). This condition is one of the main approach in establishing PAC-Bayesian results, and has been used in studies such as \cite{alquier2011PAC}, \cite{alquier2013sparse}, and \cite{mai2023bilinear}.

Following \cite{kohler2021rate,langer2021analysis} by taking the architecture of the neural networks as 
\begin{equation}
	\label{eq_langer_architecture}
	L \asymp \log (n) \,
	,
	D \asymp   n^{\frac{d}{2(2\beta+d)}}  
	,
\end{equation}
we obtain the following results.

\begin{theorem}
	\label{cor_regression_1}
	Assume that Assumption \ref{asm1} and \ref{assu_bounded} are satisfied. Put  $ C^*_{\sigma, \varsigma} := 16[C^2+\sigma^2+ C(\varsigma \lor(2C) ] $ and $ \lambda = n/ C^*_{\sigma, \varsigma} $.
	We consider the architectures as in \eqref{eq_langer_architecture}, then
	\begin{equation*}
		\mathbb{E} \, \mathbb{E}_{\theta\sim\hat{\rho}_{\lambda}} [R(\theta) ]- R(f_0)
		\leq 
		\min_{\theta\in \Theta_{L,D} }
		\frac{8}{3} [R(\theta) - R(f_0)] 
		+ 
		\mathcal{C } \psi_n
		,
	\end{equation*}
	and with any $ \delta \in (0,1) $ that:
	\begin{equation*}
		\mathbb{P}\Big\{\!  
		\mathbb{E}_{\theta\sim\hat{\rho}_{\lambda}} [R(\theta) ]- R(f_0)
		\leq
		\min_{\theta\in \Theta_{L,D} }
		4 [R(\theta) - R(f_0)] 
		+ 
		\mathcal{C }  \psi_n
		+ 	
		\mathcal{C }
		\frac{ \log(2/\delta)  }{n} 
		\! \Big\}
		\geq
		1-\delta
		,
	\end{equation*}
	where $ \psi_n := n^{-\frac{-2\beta}{2\beta +d}} 
	(\log n)^3 $ and $ \mathcal{C } >0 $ is some universal constants depending only on $ B,C, \sigma, \varsigma $.
\end{theorem}

\begin{remark}
	A similar rate is obtained in \cite{kohler2021rate} for fully connected DNNs but with only ReLU activation function. The results is then extended to sigmoid activation function in \cite{langer2021analysis}. So our work extend these works by considering a large class of activation functions that being 1-Lipschitz.	
\end{remark}

Assuming that the true underlying function \( f_0 \) belongs to a Besov space and there exists a \( \theta^* \) such that \( R(f_0) = R(\theta^*) = \min_{\theta \in \Theta_{L,D}} R(\theta) \), we can directly derive the following corollary.

\begin{corollary}
	\label{cor_regression_1_cor}
	Assume that Theorem \ref{cor_regression_1} hold true and $ \psi_n := n^{-\frac{-2\beta}{2\beta +d}} 
	(\log n)^3 $, then
	$$ \,
	\mathbb{E} \, \mathbb{E}_{\theta\sim\hat{\rho}_{\lambda}} [R(\theta) ]- R(f_0)
	\leq 
	\mathcal{C } \psi_n
	, \quad
	$$
	$$
	\mathbb{P}  \big\{\!  
	\mathbb{E}_{\theta\sim\hat{\rho}_{\lambda}} [R(\theta) ]- R(f_0)
	\leq
	\mathcal{C } 
	\psi_n
	+ 	
	\mathcal{C }
	\frac{ \log(2/\delta)  }{n} 
	\! \big\}
	\geq         1-\delta
	,
	$$
	with any $ \delta \in (0,1) $ and $ \mathcal{C } >0 $ is some universal constants depending only on $ B,C, \sigma, \varsigma $.
\end{corollary}

\begin{remark}		
	The rate $ n^{\frac{ - 2\beta}{2\beta +d}}  $ up to a log-term is similar to \cite{suzuki2019adaptivity,lee2022asymptotic} for which the authors consider sparse DNNs in Besov space. This rate is known to be minimax-optimal in the Besov spaces \citep{kerkyacharian1992density,donoho1998minimax}. However, it should be noted that both the references \cite{suzuki2019adaptivity} and \cite{lee2022asymptotic} only consider ReLU activation function.
\end{remark}

\begin{remark}	
	Our results extend the work of \cite{tinsi2022risk}, which is limited to neural networks with a single hidden layer. In contrast, we demonstrate that similar results hold for deep neural networks. We also emphasize that the methodology in \cite{tinsi2022risk} differs from ours, as their approach relies on a mirror averaging aggregation procedure (see also \cite{dalalyan2012mirror}).
\end{remark}

\begin{remark}
	With respect to the additional logarithmic factor, \cite{suzuki2019adaptivity} achieved a bound involving \( \log^2(n) \), which is slightly tighter than our current \( \log^3(n) \) rate. It is worth mentioning that \cite{lee2022asymptotic} also reported a \( \log^3(n) \) term. 
	
	We highlight that, by employing the following network architecture:
	\begin{equation}
		\label{eq_tienmt_architecture}
		L \asymp [n/\log (n) ]^{\frac{d}{4(2\beta+d)}} \,, \quad
		D \asymp [n/\log (n) ]^{\frac{d}{4(2\beta+d)}},
	\end{equation}
	we are able to significantly reduce the logarithmic factor to \( \left[
	\log (n) 
	\right]^{\frac{2\beta}{2\beta+d} } \), see Proposition \ref{cor_regression_2} below. To the best of our knowledge, this represents a notable improvement.
\end{remark}

\begin{prop}
	\label{cor_regression_2}
	Assume the same conditions for Theorem \ref{cor_regression_1} hold.
	We consider the architectures as in \eqref{eq_tienmt_architecture}, then
	\begin{equation*}
		\mathbb{E} \, \mathbb{E}_{\theta\sim\hat{\rho}_{\lambda}} [R(\theta) ]- R(f_0)
		\leq 
		\min_{\theta\in \Theta_{L,D} }
		\frac{8}{3} [R(\theta) - R(f_0)] 
		+ 
		\mathcal{C } \tilde{\psi}_n
		,
	\end{equation*}
	and with any $ \delta \in (0,1) $ that:
	\begin{equation*}
		\mathbb{P}\Biggl\{\!  
		\mathbb{E}_{\theta\sim\hat{\rho}_{\lambda}} [R(\theta) ]- R(f_0)
		\leq
		\min_{\theta\in \Theta_{L,D} }
		\left[
		4 [R(\theta) - R(f_0)] 
		+ 
		\mathcal{C } 
		\tilde{\psi}_n
		+ 	
		\mathcal{C }
		\frac{ \log(2/\delta)  }{n} 
		\right]
		\! \Biggr\}
		\geq
		1-\delta
		,
	\end{equation*}
	where $ \tilde{\psi}_n := 		\left[
	\log (n) / n  
	\right]^{\frac{2\beta}{2\beta+d} } $ and $ \mathcal{C } >0 $ is some universal constants depending only on $ B,C, \sigma, \varsigma $.
\end{prop}

\section{Results for binary classification}
\label{sc_classification}
\subsection{Setting}
Now,  a binary logistic regression  model is considered. More specially, let $ (Y_i,X_i)\in \lbrace -1,1 \rbrace \times [0,1]^d $, $ i=1,\ldots,n $, be a collection of $n$ \text{i.i.d.} random samples drawn from a joint distribution $ P $ satisfying for all $x \in [0,1]^d  $,
\begin{equation*}
	Y_i \mid X_i=x~\sim 2\mathcal{B}\big( \eta(x) \big) - 1, ~ \text{ with } ~ \eta(x) = P(Y_i = 1| X_i = x)
	, 
\end{equation*}
where $\mathcal{B}( p )$ denotes the Bernoulli distribution with parameter $ p $. Here $ \eta(x) = (1+ e^{-f_0 (X)} )^{-1} $ and $f_0:[0,1]^d\rightarrow \mathbb{R}$ is the true unknown function.

We focus on the cross-entropy (logistic) loss function given by,
$$
\ell(y, u) = \log(1 + e^{-y u}) 
$$ for all $(y,u) \in \{-1, 1\} \times \mathbb{R} $.
For any predictor $ f_\theta $, we define the prediction risk as,
\[ R(\theta) = \mathbb{E} \big[ \ell \big(Y_1, f_\theta (X_1)  \big)  \big] 
. 
\]
and its empirical risk as $ r_{n}( \theta )
=
\frac{1}{n}\sum_{i=1}^{n}
\ell (Y_{i} , f_\theta (X_{i}) ) $.

In a similar manner, we investigate the exponentially weighted aggregate (EWA) procedure as a potential solution to this problem. This involves examining the Gibbs posterior distribution, a key tool for deriving probabilistic model weights based on their associated risks or losses, defined in \eqref{eq_mainporsterior}, $ \hat{\rho}_\lambda(\theta) $.

\subsection{Excess risk bounds for cross-entropy loss}

Let $\theta^*$ denote a minimizer of $R$ when it exists:
$
R(\theta^*) 
= 
\min_{\theta\in \Theta_{L,D} } R(\theta). 
$. We first formulate an assumption for the context of logistic regression.
\begin{assume}
	\label{assum_bernstein}
	Assume that for any $\theta\in\Theta $, there is a constant $K>0 $ such that, 
	$$
	\mathbb{E} \{ 
	[ \ell (Y_i, f_\theta(X_i)) - \ell (Y_i, f_{0} (X_i)) ]^2
	\}
	\leq 
	K [ R (\theta) - R(f_0) ]
	.
	$$
\end{assume}

Assumption \ref{assum_bernstein} is known as a kind of Bernstein's type condition in learning theory, see e.g. \cite{bartlett2006convexity} and \cite{alquier2021user} for some discussion. This assumption is used to obtain a fast learning rate.

Now, we state our main theorem for fully connected DNN in the context of binary logistic regression using cross-entropy loss.

\begin{theorem}
	\label{cor_classification_cor}
	We consider the architecture for DNN as  $
	L \asymp [n/\log (n) ]^{\frac{d}{4(\beta+d)}} 
	\,, 
	D \asymp [n/\log (n) ]^{\frac{d}{4(\beta+d)}}
	$.
	Assume that Assumption \ref{asm1} and \ref{assum_bernstein} are satisfied. Take $\lambda= n/ \max(2K,C) $, we have:
	\begin{equation*}
		\mathbb{E} \mathbb{E}_{\theta\sim\hat{\rho}_{\lambda}} [R(\theta) ]- R(f_0)
		\leq 
		\min_{\theta\in \Theta_{L,D} }
		2 [R(\theta) - R(f_0)] 
		+
		\mathcal{H }\, \psi_n
		,
	\end{equation*}
	and with any $ \delta \in (0,1) $ that:
	\begin{equation*}
		\mathbb{P}\Biggl\{\!  
		\mathbb{E}_{\theta\sim\hat{\rho}_{\lambda}} [R(\theta) ]- R(f_0)
		\leq
		\min_{\theta\in \Theta_{L,D} }
		3 [R(\theta) - R(f_0)] 
		+
		\mathcal{H } \psi_n
		+ 	
		\mathcal{H }
		\frac{ \log(2/\delta)  }{n} 
		\! \Biggr\}
		\geq
		1-\delta
		,
	\end{equation*}
	where 
	$ \psi_n := 		\left[
	\frac{  \log (n)}{n}  
	\right]^{\frac{\beta}{\beta+d} }
	$ 
	and $ \mathcal{H } >0 $ is some universal constants depending only on $ B,C, K $.
\end{theorem}

\begin{remark}
	The excess error rate $\psi_n = \left[ \log(n)/n \right]^{\frac{\beta}{\beta + d}}$, up to a logarithmic factor, is known to be minimax-optimal for the logistic (cross-entropy) loss over a class of Hölder-smooth functions. This is justified by the minimax lower bound established in Corollary 2.1 of \cite{zhang2024classification}. 
	\\
	A corresponding upper bound of order $\left[ \log^5(n)/n \right]^{\frac{\beta}{\beta + d}}$ was also derived in Theorem 2.2 (equation 2.17) of \cite{zhang2024classification}, for sparse deep neural networks with ReLU activation functions. 
	Therefore, our results generalize these findings to a broader setting involving fully connected neural networks with general activation functions. 
	Notably, we also achieve a substantial improvement by reducing the logarithmic term in the convergence rate.
\end{remark}

If we assume that the true underlying function $ f_0 $ is in Besov space and that there exist a $ \theta^* $ such that $ R(f_0) = R(\theta^*) = \min_{\theta\in \Theta_{L,D} } R(\theta) $, then we immediately obtain the following corollary from Theorem  \ref{cor_classification_cor}.

\begin{corollary}
	\label{cor_classification}
	Assume that the true underlying function $f_0$ belongs to a Besov space. and that there exists $\theta^* \in \Theta_{L,D}$ such that $R(f_0) = R(\theta^*) = \min_{\theta \in \Theta_{L,D}} R(\theta)$. From Theorem  \ref{cor_classification_cor}, we obtain the following result:
	$$
	\mathbb{E} \,
	\mathbb{E}_{\theta\sim\hat{\rho}_{\lambda}} [R(\theta) ]- R(f_0)
	\leq 
	\mathcal{H }\, \psi_n
	,
	$$
	$$
	\mathbb{P}\Bigg\{\!  
	\mathbb{E}_{\theta\sim\hat{\rho}_{\lambda}} [R(\theta) ]- R(f_0)
	\leq
	\mathcal{H } \psi_n
	+ 	
	\mathcal{H }	\frac{ \log(2/\delta)  }{n} 
	\! \Bigg\}
	\geq  	1-\delta
	,
	$$
	with any $ \delta \in (0,1) $  and $ \mathcal{H } >0 $ is some universal constants depending only on $ B,C, K $.
\end{corollary}

As mentioned in the introduction, most theoretical work on Bayesian deep neural networks focuses either on the nonparametric regression setup or on the posterior contraction rate in logistic regression. Our results are novel in this setting as they focus on the prediction risk of probabilistic DNNs.

\subsection{Misclassification excess risk bounds}

Up to now, we have just worked with the logistic/cross entropy loss meaning that the excess risk is only for this loss. However, in classification, the misclassification risk is also of very important to quantify the accuracy of a classifier. We tackle that problem in this section.

In the logistic model, we define our classifier $ \hat{\eta}^{DNN}_\theta (x) $ as
\begin{equation*}
	\begin{cases}
		\hat{\eta}^{DNN}_\theta (X) = 1 , & \text{ if }   (1+ e^{-f_\theta (X)} )^{-1} \geq 1/2 
		,
		\\
		\hat{\eta}^{DNN}_\theta (X) = -1 , & \text{ if }    
		(1+ e^{-f_\theta (X)} )^{-1} < 1/2 
		.
	\end{cases}
\end{equation*}
More specifically, it is noted that $  (1+ e^{-f_\theta (X)} )^{-1} \geq 1/2  $ is equivalent to $ f_\theta (X) \geq 0 $. Thus, the classifier $ \hat{\eta}^{DNN}_\theta (x) $ reduces to as $ {\rm sign} (f_\theta (X)) $.

The accuracy of a classifier \( \hat{\eta}^{DNN}_\theta (x) \) is determined by the prediction or misclassification error, defined as
$$
R_{0/1} (\theta) = \mathbb{P}
[ Y \neq \hat{\eta}^{DNN}_\theta (x)  ]
.
$$
Let $\theta^*$ denote a minimizer of $R_{0/1} $ when it exists:
$
R_{0/1}^* 
= 
\min_{\theta\in \Theta_{L,D} } R_{0/1} (\theta). 
$
and consider its misclassification excess risk 
\[
\mathbb{E} R_{0/1}(\theta) - R^*_{0/1}
.
\]

Here after, we introduce the following condition on the underlying classifier.

\begin{assume}
	\label{assum_margin}
	Assume that there exist $ C_{m.g} >0 $ such that 
	$
	\mathbb{P} ( |\eta(x) - 1/2| \leq h ) \leq C_{m.g} h
	,
	$
	for all $ 0 <h<1/2 $.
\end{assume}

Assumption \ref{assum_margin} introduces a low-noise condition commonly employed in the classification literature, as in \citep{abramovich2018high,tsybakov2004optimal,mammen1999smooth,bartlett2006convexity}. Near the decision boundary \( \{x : \eta(x) = 1/2\} \), which in the context of logistic regression corresponds to the hyperplane \( f(x) = 0 \) where \( \eta(x) = (1 + e^{-f(x)})^{-1} \), classifiers have difficulty in returning a good prediction. Predicting the class label in this region is challenging due to the predominance of noise in the label information. As a result, it is a reasonable assumption that \( \eta(x) \) does not frequently take values close to \( 1/2 \).

\begin{theorem}
	\label{thm_mis_claass_excess}
	Assume that Corollary \ref{cor_classification} holds and that Assumption \ref{assum_margin} holds. We consider the architecture as  
	$	 L \asymp [n/\log (n) ]^{\frac{\beta+d}{4(2\beta+d)}} 
	\,, \,
	D \asymp [n/\log (n) ]^{\frac{\beta+d}{4(2\beta+d)}}
	$, 
	then
	\begin{equation*}
		\mathbb{E} \mathbb{E}_{\theta\sim\hat{\rho}_{\lambda}} [R_{0/1} (\theta) ]- R_{0/1}(f_0)
		\leq 
		\mathcal{C } \psi_n
		.
	\end{equation*}
	and with any $ \delta \in (0,1) $ that:
	\begin{equation*}
		\mathbb{P}\Biggl\{\!  
		\mathbb{E}_{\theta\sim\hat{\rho}_{\lambda}} [R_{0/1} (\theta) ]- R_{0/1} (f_0)
		\leq
		\mathcal{C } \psi_n
		+ 	
		\mathcal{C }
		\frac{ \log(2/\delta)  }{n} 
		\! \Biggr\}
		\geq
		1-\delta
		,
	\end{equation*}
	with $ \psi_n 
	:= \left[
	\log (n) / n  
	\right]^{\frac{ \beta}{2\beta+d} } $, where $ \mathcal{C } >0 $ is some universal constants depending only on $ B,C, K,C_{m.g} $.
\end{theorem}

\begin{remark}
	In \cite{zhang2024classification}, the authors derived a slower misclassification excess risk of order $ \left[ \log^5(n)/n \right]^{\frac{\beta}{2\beta + 2d}}$ under the logistic loss, (see in Theorem 2.2 in \cite{zhang2024classification}). 
	Under the standard Hölder smoothness condition, the minimax-optimal rate for misclassification excess risk is known to be of order $ n^{-\beta/(2\beta + d)}$, see \cite{audibert2007fast}. 
	Thus, our result achieves the minimax-optimal rate up to a logarithmic factor. 
	Minima-optimal rate for misclassification risk has also been obtained in previous work using hinge loss in \cite{kim2021fast} for sparse DNNs with ReLU activation function.
\end{remark}

\begin{remark}
	Furthermore, the logarithmic term in our bound improves upon the result in \cite{kim2021fast}, where a rate of $  \left[ \log^3(n)/n \right]^{\frac{\beta}{2\beta + d}}$ was obtained under the hinge loss and with ReLU activation functions only (see Theorem 3.3 in \cite{kim2021fast}).
\end{remark}

\begin{remark}
	In \cite{mai2024misclassification}, which employs the hinge loss and focuses on sparse DNNs, the author derived a similar convergence rate of order $n^{-\frac{\beta}{2\beta + d}} \log(n)$ for a strictly noiseless setting, but restricted attention to classifiers within the $\beta$-Hölder function class. 
	However, \cite{mai2024misclassification} did not establish a sharp oracle inequality as we do in this work. 
	Therefore, our results represent a notable contribution by providing a minimax-optimal sharp oracle inequality. 
	A sharp oracle inequality not only bounds the excess risk but also matches the leading constant under certain conditions, giving a more precise characterization of the estimator’s performance and providing stronger theoretical guarantees. 
	Our log-term $	 \log^{\frac{\beta}{2\beta+d} } (n) $ also slightly improves over \cite{mai2024misclassification}. 
\end{remark}

\section{Conclusion}

In this paper, we established PAC-Bayesian generalization bounds for fully connected Bayesian deep neural networks under a Gaussian prior, a prior choice widely adopted in both theoretical and applied contexts. Focusing on an exponentially weighted aggregate estimator, we derived upper bounds on the prediction risk in both nonparametric regression and binary classification with logistic loss. These bounds achieve, up to logarithmic factors, the minimax-optimal convergence rates, thereby providing strong theoretical guarantees on the performance of deep neural network probabilistic models.

A key novelty of our work lies in its focus on fully connected architectures, as opposed to the predominant emphasis on sparse or shallow networks in the existing PAC-Bayesian literature. 
Furthermore, we contribute  sharp oracle inequalities for classification—something that previous related works, such as \cite{mai2024misclassification}, did not attain—thus offering a more refined characterization of the excess risk and narrowing the gap between theoretical upper bounds and practical performance. Moreover, the additional logarithmic term is showed to be better than existing results.

Our results not only extend prior work limited to regression or classification alone, but also unify the treatment of both learning tasks under a single framework. We hope that these findings will stimulate further theoretical investigations into different architectures, alternative prior structures, and extensions to other learning settings, such as multi-class classification or structured prediction under uncertainty.

\subsubsection*{Acknowledgments}
The views, results, and opinions expressed in this work are solely those of the author and do not, in any way, represent those of the Norwegian Institute of Public Health.

\subsubsection*{Conflicts of interest/Competing interests}
The author declares no potential conflict of interests.

\appendix
\section{Proof}
\label{sc_proof}

\subsection{Auxiliary lemmas}

\begin{lemma}[Bernstein's inequality, Theorem 5.2.1 in~\cite{catoni2004statistical}]
	\label{lemma:bernstein}
Define the function $g$ as $g(0) = 1$ and, for $x\neq 0$,
	$ g(x) = \frac{{\rm e}^x - 1 - x}{x^2}. $
	Let $U_1,\dots,U_n$ be i.i.d random variables such that $\mathbb{E}(U_i)$ is well defined and $U_i- \mathbb{E}(U_i) \leq C$ almost surely for some $C\in\mathbb{R}$. Then
	$$
	\mathbb{E}\left( {\rm e}^{t \sum_{i=1}^n [U_i - \mathbb{E}(U_i) ] } \right)
	\leq {\rm e}^{ g\left(C t\right) n t^2 {\rm Var}(U_i) }.
	$$
\end{lemma}

\begin{lemma}[Donsker and Varadhan's variational formula, \cite{catonibook}]
	\label{lemma:dv}
	For any measurable, bounded function $h:\Theta\rightarrow\mathbb{R}$ we have:
	\begin{equation*}
	\log \mathbb{E}_{\theta\sim\pi}\left[{\rm e}^{h(\theta)} \right] 
    =
\sup_{\nu \in\mathcal{P}(\Theta)}
\Bigl[\mathbb{E}_{\theta\sim \nu}[h(\theta)] 
-	
\mathcal{K} (\nu , \pi)\Bigr].
	\end{equation*}
	Moreover, the supremum with respect to $\rho$ in the right-hand side is
	reached for the Gibbs measure
	$\pi_{h}$ defined by its density with respect to $\pi$
	\begin{equation*}
	\frac{{\rm d}\pi_{h}}{{\rm d}\pi}(\theta) =  \frac{{\rm e}^{h(\theta)}}
	{ \mathbb{E}_{\vartheta\sim\pi}\left[{\rm e}^{h(\vartheta)} \right] }.
	\end{equation*}
\end{lemma}

We will use a version of Bernstein's inequality taken from \cite{MR2319879} (Inequality 2.21, page 24).

\begin{lemma}
	\label{lemmemassart} 
	Let $ U_{1}, \ldots,  U_{n} $ be independent real
	valued random variables and assume that there are two constants
	$v$ and $w$ such that
	$
	\sum_{i=1}^{n} \mathbb{E}[ U_{i}^{2}] \leq v
	,
	$
	and for all integers $k\geq3$,
	$
	\sum_{i=1}^{n} \mathbb{E}\left[( U_{i})^{k}_+ \right] \leq v\frac
	{k!w^{k-2}}{2}.
	$
	Then, for any $\zeta\in(0,1/w)$,
	\[
	\mathbb{E}
	\exp\left[\zeta\sum_{i=1}^{n}\left[ U_{i} - \mathbb{E}( U_{i})\right]
	\right]
	\leq\exp\left(\frac{v\zeta^{2}}{2(1-w\zeta)} \right) .
	\]
\end{lemma}


\begin{thebibliography}{}
	
	\bibitem[Abramovich and Grinshtein, 2018]{abramovich2018high}
	Abramovich, F. and Grinshtein, V. (2018).
	\newblock High-dimensional classification by sparse logistic regression.
	\newblock {\em IEEE Transactions on Information Theory}, 65(5):3068--3079.
	
	\bibitem[Alquier, 2024]{alquier2021user}
	Alquier, P. (2024).
	\newblock User-friendly introduction to {PAC-B}ayes bounds.
	\newblock {\em Foundations and Trends{\textregistered} in Machine Learning},
	17(2):174--303.
	
	\bibitem[Alquier and Biau, 2013]{alquier2013sparse}
	Alquier, P. and Biau, G. (2013).
	\newblock Sparse single-index model.
	\newblock {\em J. Mach. Learn. Res.}, 14:243--280.
	
	\bibitem[Alquier and Kengne, 2024]{alquier2024minimaxopti}
	Alquier, P. and Kengne, W. (2024).
	\newblock {Minimax optimality of deep neural networks on dependent data via
		PAC-Bayes bounds}.
	\newblock {\em arXiv:2410.21702}.
	
	\bibitem[Alquier and Lounici, 2011]{alquier2011PAC}
	Alquier, P. and Lounici, K. (2011).
	\newblock P{AC}-{B}ayesian bounds for sparse regression estimation with
	exponential weights.
	\newblock {\em Electron. J. Stat.}, 5:127--145.
	
	\bibitem[Anthony and Bartlett, 2009]{anthony1999neural}
	Anthony, M. and Bartlett, P.~L. (2009).
	\newblock {\em Neural network learning: Theoretical foundations}, volume~9.
	\newblock cambridge university press, Cambridge.
	
	\bibitem[Bai et~al., 2020]{bai2020efficient}
	Bai, J., Song, Q., and Cheng, G. (2020).
	\newblock Efficient variational inference for sparse deep learning with
	theoretical guarantee.
	\newblock {\em Advances in Neural Information Processing Systems}, 33:466--476.
	
	\bibitem[Barron, 1994]{Barron94estimation}
	Barron, A.~R. (1994).
	\newblock Approximation and estimation bounds for artificial neural networks.
	\newblock {\em Machine learning}, 14:115--133.
	
	\bibitem[Bartlett et~al., 2017]{Bartlett2017MarginBoundsNNs}
	Bartlett, P.~L., Foster, D.~J., and Telgarsky, M.~J. (2017).
	\newblock Spectrally-normalized margin bounds for neural networks.
	\newblock {\em Advances in neural information processing systems}, 30.
	
	\bibitem[Bartlett et~al., 2006]{bartlett2006convexity}
	Bartlett, P.~L., Jordan, M.~I., and McAuliffe, J.~D. (2006).
	\newblock Convexity, classification, and risk bounds.
	\newblock {\em Journal of the American Statistical Association},
	101(473):138--156.
	
	\bibitem[Bhattacharya et~al., 2024]{bhattacharya2024comprehensive}
	Bhattacharya, S., Liu, Z., and Maiti, T. (2024).
	\newblock Comprehensive study of variational bayes classification for dense
	deep neural networks.
	\newblock {\em Statistics and Computing}, 34(1):17.
	
	\bibitem[Bingham et~al., 2019]{bingham2019pyro}
	Bingham, E., Chen, J.~P., Jankowiak, M., Obermeyer, F., Pradhan, N.,
	Karaletsos, T., Singh, R., Szerlip, P., Horsfall, P., and Goodman, N.~D.
	(2019).
	\newblock Pyro: Deep universal probabilistic programming.
	\newblock {\em Journal of machine learning research}, 20(28):1--6.
	
	\bibitem[Boucheron et~al., 2013]{boucheron2013concentration}
	Boucheron, S., Lugosi, G., and Massart, P. (2013).
	\newblock {\em Concentration inequalities: A nonasymptotic theory of
		independence}.
	\newblock Oxford University Press, Oxford.
	
	\bibitem[Catoni, 2004]{catoni2004statistical}
	Catoni, O. (2004).
	\newblock {\em Statistical learning theory and stochastic optimization}, volume
	1851 of {\em Saint-Flour Summer School on Probability Theory 2001 (Jean
		Picard ed.), Lecture Notes in Mathematics}.
	\newblock Springer-Verlag, Berlin.
	
	\bibitem[Catoni, 2007]{catonibook}
	Catoni, O. (2007).
	\newblock {\em {PAC}-{B}ayesian supervised classification: the thermodynamics
		of statistical learning}.
	\newblock IMS Lecture Notes---Monograph Series, 56. Institute of Mathematical
	Statistics, Beachwood, OH.
	
	\bibitem[Ch{\'e}rief-Abdellatif, 2020]{cherief2020convergence}
	Ch{\'e}rief-Abdellatif, B.-E. (2020).
	\newblock Convergence rates of variational inference in sparse deep learning.
	\newblock In III, H.~D. and Singh, A., editors, {\em Proceedings of the 37th
		International Conference on Machine Learning}, volume 119, pages 1831--1842.
	PMLR.
	
	\bibitem[Dalalyan and Tsybakov, 2012]{dalalyan2012mirror}
	Dalalyan, A.~S. and Tsybakov, A. (2012).
	\newblock Mirror averaging with sparsity priors.
	\newblock {\em Bernoulli}, 18(3):914--944.
	
	\bibitem[Donoho and Johnstone, 1998]{donoho1998minimax}
	Donoho, D.~L. and Johnstone, I.~M. (1998).
	\newblock Minimax estimation via wavelet shrinkage.
	\newblock {\em The annals of Statistics}, 26(3):879--921.
	
	\bibitem[Elbr{\"a}chter et~al., 2021]{elbrachter2021deep}
	Elbr{\"a}chter, D., Perekrestenko, D., Grohs, P., and B{\"o}lcskei, H. (2021).
	\newblock Deep neural network approximation theory.
	\newblock {\em IEEE Transactions on Information Theory}, 67(5):2581--2623.
	
	\bibitem[Gin{\'e} and Nickl, 2016]{gine2016mathematical}
	Gin{\'e}, E. and Nickl, R. (2016).
	\newblock {\em Mathematical foundations of infinite-dimensional statistical
		models}, volume~40.
	\newblock Cambridge university press.
	
	\bibitem[Goodfellow et~al., 2016]{goodfellow2016deep}
	Goodfellow, I., Bengio, Y., and Courville, A. (2016).
	\newblock {\em Deep learning}.
	\newblock MIT press.
	
	\bibitem[Guedj, 2019]{guedj2019primer}
	Guedj, B. (2019).
	\newblock A primer on {PAC}-{B}ayesian learning.
	\newblock In {\em S{MF} 2018: {C}ongr\`es de la {S}oci\'{e}t\'{e}
		{M}ath\'{e}matique de {F}rance}, volume~33 of {\em S\'{e}min. Congr.}, pages
	391--413. Soc. Math. France.
	
	\bibitem[Hayakawa and Suzuki, 2020]{Suzuki2019Superiority}
	Hayakawa, S. and Suzuki, T. (2020).
	\newblock On the minimax optimality and superiority of deep neural network
	learning over sparse parameter spaces.
	\newblock {\em Neural Networks}, 123:343--361.
	
	\bibitem[Haykin, 1998]{haykin1998neural}
	Haykin, S. (1998).
	\newblock {\em Neural networks: a comprehensive foundation}.
	\newblock Prentice Hall PTR.
	
	\bibitem[Hertz, 2018]{hertz2018introduction}
	Hertz, J.~A. (2018).
	\newblock {\em Introduction to the theory of neural computation}.
	\newblock Crc Press.
	
	\bibitem[Imaizumi and Fukumizu, 2019]{Imaizumi19DNN}
	Imaizumi, M. and Fukumizu, K. (2019).
	\newblock Deep neural networks learn non-smooth functions effectively.
	\newblock In {\em The 22nd international conference on artificial intelligence
		and statistics}, pages 869--878. PMLR.
	
	\bibitem[Jantre et~al., 2023]{jantre2023layer}
	Jantre, S., Bhattacharya, S., and Maiti, T. (2023).
	\newblock Layer adaptive node selection in bayesian neural networks:
	Statistical guarantees and implementation details.
	\newblock {\em Neural Networks}, 167:309--330.
	
	\bibitem[Jantre et~al., 2024]{jantre2024spike}
	Jantre, S., Bhattacharya, S., and Maiti, T. (2024).
	\newblock Spike-and-slab shrinkage priors for structurally sparse Bayesian neural networks
	\newblock {\em IEEE Transactions on Neural Networks and Learning Systems}, 2024.
	
	\bibitem[Kerkyacharian and Picard, 1992]{kerkyacharian1992density}
	Kerkyacharian, G. and Picard, D. (1992).
	\newblock Density estimation in besov spaces.
	\newblock {\em Statistics \& probability letters}, 13:15--24.
	
	\bibitem[Kim et~al., 2021]{kim2021fast}
	Kim, Y., Ohn, I., and Kim, D. (2021).
	\newblock Fast convergence rates of deep neural networks for classification.
	\newblock {\em Neural Networks}, 138:179--197.
	
	\bibitem[Kohler and Langer, 2021]{kohler2021rate}
	Kohler, M. and Langer, S. (2021).
	\newblock On the rate of convergence of fully connected deep neural network
	regression estimates.
	\newblock {\em The Annals of Statistics}, 49(4):2231--2249.
	
	\bibitem[Kohler et~al., 2023]{kohler2023estimation}
	Kohler, M., Langer, S., and Reif, U. (2023).
	\newblock Estimation of a regression function on a manifold by fully connected
	deep neural networks.
	\newblock {\em Journal of Statistical Planning and Inference}, 222:160--181.
	
	\bibitem[Kong and Kim, 2024]{kong2024posterior}
	Kong, I. and Kim, Y. (2024).
	\newblock Posterior concentrations of fully-connected bayesian neural networks
	with general priors on the weights.
	\newblock {\em arXiv preprint arXiv:2403.14225}.
	
	\bibitem[Kong et~al., 2023]{kong2023masked}
	Kong, I., Yang, D., Lee, J., Ohn, I., Baek, G., and Kim, Y. (2023).
	\newblock Masked bayesian neural networks: Theoretical guarantee and its
	posterior inference.
	\newblock In {\em International Conference on Machine Learning}, pages
	17462--17491. PMLR.
	
	\bibitem[Langer, 2021]{langer2021analysis}
	Langer, S. (2021).
	\newblock Analysis of the rate of convergence of fully connected deep neural
	network regression estimates with smooth activation function.
	\newblock {\em Journal of Multivariate Analysis}, 182:104695.
	
	\bibitem[LeCun et~al., 2015]{lecun2015deep}
	LeCun, Y., Bengio, Y., and Hinton, G. (2015).
	\newblock Deep learning.
	\newblock {\em Nature}, 521(7553):436--444.
	
	\bibitem[Lee and Lee, 2022]{lee2022asymptotic}
	Lee, K. and Lee, J. (2022).
	\newblock {Asymptotic properties for Bayesian neural network in Besov space}.
	\newblock {\em Advances in Neural Information Processing Systems},
	35:5641--5653.
	
	\bibitem[Liu, 2021]{liu2021variable}
	Liu, J. (2021).
	\newblock Variable selection with rigorous uncertainty quantification using
	deep bayesian neural networks: Posterior concentration and bernstein-von
	mises phenomenon.
	\newblock In {\em International Conference on Artificial Intelligence and
		Statistics}, pages 3124--3132. PMLR.
	
	\bibitem[Mai, 2023a]{mai2023bilinear}
	Mai, T.~T. (2023a).
	\newblock From bilinear regression to inductive matrix completion: a
	quasi-bayesian analysis.
	\newblock {\em Entropy}, 25(2):333.
	
	\bibitem[Mai, 2023b]{mai2023reduced}
	Mai, T.~T. (2023b).
	\newblock A reduced-rank approach to predicting multiple binary responses
	through machine learning.
	\newblock {\em Statistics and Computing}, 33(6):136.
	
	\bibitem[Mai, 2024a]{mai2024sparse}
	Mai, T.~T. (2024a).
	\newblock {A sparse PAC-Bayesian approach for high-dimensional quantile
		prediction}.
	\newblock {\em arXiv preprint arXiv:2409.01687}.
	
	\bibitem[Mai, 2024b]{mai2024high}
	Mai, T.~T. (2024b).
	\newblock High-dimensional sparse classification using exponential weighting
	with empirical hinge loss.
	\newblock {\em Statistica Neerlandica}, 78(4):664--691.
	
	\bibitem[Mai, 2025]{mai2024misclassification}
	Mai, T.~T. (2025).
	\newblock {Misclassification bounds for PAC-Bayesian sparse deep learning}.
	\newblock {\em Machine Learning}, 114(1):18.
	
	\bibitem[Mammen and Tsybakov, 1999]{mammen1999smooth}
	Mammen, E. and Tsybakov, A.~B. (1999).
	\newblock Smooth discrimination analysis.
	\newblock {\em The Annals of Statistics}, 27(6):1808--1829.
	
	\bibitem[Massart, 2007]{MR2319879}
	Massart, P. (2007).
	\newblock {\em Concentration inequalities and model selection}, volume 1896 of
	{\em Lecture Notes in Mathematics}.
	\newblock Springer, Berlin.
	\newblock Lectures from the 33rd Summer School on Probability Theory held in
	Saint-Flour, July 6--23, 2003, Edited by Jean Picard.
	
	\bibitem[McAllester, 1998]{McA}
	McAllester, D. (1998).
	\newblock Some {{PAC}}-{B}ayesian theorems.
	\newblock In {\em Proceedings of the Eleventh Annual Conference on
		Computational Learning Theory}, pages 230--234, New York. ACM.
	
	\bibitem[Neyshabur et~al., 2018]{Srebro2018PACMarginBoundsNNs}
	Neyshabur, B., Bhojanapalli, S., and Srebro, N. (2018).
	\newblock A pac-bayesian approach to spectrally-normalized margin bounds for
	neural networks.
	\newblock In {\em International Conference on Learning Representations}.
	
	\bibitem[Polson and Ro\v{c}kov\'{a}, 2018]{Rockova2018}
	Polson, N.~G. and Ro\v{c}kov\'{a}, V. (2018).
	\newblock Posterior concentration for sparse deep learning.
	\newblock In {\em Advances in Neural Information Processing Systems},
	volume~31, Montreeal, Canada. Curran Associates, Inc.
	
	\bibitem[Ripley, 2007]{ripley2007pattern}
	Ripley, B.~D. (2007).
	\newblock {\em Pattern recognition and neural networks}.
	\newblock Cambridge university press.
	
	\bibitem[Sawano, 2018]{sawano2018theory}
	Sawano, Y. (2018).
	\newblock {\em Theory of Besov spaces}, volume~56.
	\newblock Springer.
	
	\bibitem[Schmidt-Hieber, 2020]{SchmidtHieberDNN}
	Schmidt-Hieber, A.~J. (2020).
	\newblock Nonparametric regression using deep neural networks with relu
	activation function.
	\newblock {\em Annals of statistics}, 48(4):1875--1897.
	
	\bibitem[Shawe-Taylor and Williamson, 1997]{STW}
	Shawe-Taylor, J. and Williamson, R. (1997).
	\newblock A {{PAC}} analysis of a {B}ayes estimator.
	\newblock In {\em Proceedings of the Tenth Annual Conference on Computational
		Learning Theory}, pages 2--9, New York. ACM.
	
	\bibitem[Steffen and Trabs, 2022]{steffen2022pac}
	Steffen, M.~F. and Trabs, M. (2022).
	\newblock Pac-bayes training for neural networks: sparsity and uncertainty
	quantification.
	\newblock {\em arXiv preprint arXiv:2204.12392}.
	
	\bibitem[Sun et~al., 2022]{sun2022learning}
	Sun, Y., Song, Q., and Liang, F. (2022).
	\newblock Learning sparse deep neural networks with a spike-and-slab prior.
	\newblock {\em Statistics \& probability letters}, 180:109246.
	
	\bibitem[Suzuki, 2018]{Suzuki18DNNkerels}
	Suzuki, T. (2018).
	\newblock Fast generalization error bound of deep learning from a kernel
	perspective.
	\newblock In {\em International Conference on Artificial Intelligence and
		Statistics}, pages 1397--1406. PMLR.
	
	\bibitem[Suzuki, 2019]{suzuki2019adaptivity}
	Suzuki, T. (2019).
	\newblock Adaptivity of deep relu network for learning in besov and mixed
	smooth besov spaces: optimal rate and curse of dimensionality.
	\newblock In {\em The 7tth International Conference on Learning
		Representations}.
	
	\bibitem[Tinsi and Dalalyan, 2022]{tinsi2022risk}
	Tinsi, L. and Dalalyan, A. (2022).
	\newblock Risk bounds for aggregated shallow neural networks using gaussian
	priors.
	\newblock In {\em Proceedings of the 35th Conference on Learning Theory},
	volume PMLR:178, pages 227--253. PMLR.
	
	\bibitem[Tsybakov, 2004]{tsybakov2004optimal}
	Tsybakov, A.~B. (2004).
	\newblock Optimal aggregation of classifiers in statistical learning.
	\newblock {\em The Annals of Statistics}, 32(1):135--166.
	
	\bibitem[Vladimirova et~al., 2019]{Vladimirova2019PriorsBNNsUnits}
	Vladimirova, M., Verbeek, J., Mesejo, P., and Arbel, J. (2019).
	\newblock Understanding priors in bayesian neural networks at the unit level.
	\newblock In {\em International Conference on Machine Learning}, pages
	6458--6467. PMLR.
	
	\bibitem[Yara and Terada, 2025]{yara2024nonparametric}
	Yara, A. and Terada, Y. (2025).
	\newblock Nonparametric logistic regression with deep learning.
	\newblock {\em Bernoulli}, page in press.
	
	\bibitem[Zhang et~al., 2021]{ZhangUnderstandingDL2017}
	Zhang, C., Bengio, S., Hardt, M., Recht, B., and Vinyals, O. (2021).
	\newblock Understanding deep learning (still) requires rethinking
	generalization.
	\newblock {\em Communications of the ACM}, 64(3):107--115.
	
	\bibitem[Zhang et~al., 2024]{zhang2024classification}
	Zhang, Z., Shi, L., and Zhou, D.-X. (2024).
	\newblock Classification with deep neural networks and logistic loss.
	\newblock {\em Journal of Machine Learning Research}, 25(125):1--117.
	
\end{thebibliography}
\end{document}